\documentclass[12pt]{article}
\usepackage{amsfonts}
\usepackage{amssymb}
\usepackage{amsmath}
\usepackage[pdftex]{hyperref}

\date{}

\begin{document}

{\footnotesize \noindent
This is the author's translation  of his paper
originally published in Russian in

\emph{Diskretnyi Analiz}, issue 31,  61--70, 91 (1977),

\qquad Novosibirsk, Institute of Mathematics

\qquad of the Siberian Branch of the Academy of Sciences of the
USSR;

\url{http://www.ams.org/mathscinet-getitem?mr=543806}

\url{http://www.zentralblatt-math.org/zmath/en/advanced/?q=an:435.05025};

\noindent The Russian original can be downloaded from

\url{http://logic.pdmi.ras.ru/~yumat/papers/28_paper/28_home.html}.

\noindent
The author is very grateful to Martin Davis
for some help with the English.
}

\

\hfill UDK 519.1

\

\begin{center}
ON A CERTAIN REPRESENTATION OF THE CHROMATIC
POLYNOMIAL\\[5mm]
Yu.V.Matiyasevich
\end{center}

\

1. NOTATION. We introduce notation
for certain standard notions in graph theory.
(Definitions of some of these notions are
given below.)

By a graph we mean a non-oriented graph, possibly with multiple edges
and loops.  $n(H)$ and~$m(H)$ will denote the number of vertices and
edges, respectively, of a graph~$H$. Let $H\le G$ mean that $H$ is an
edge subgraph of a graph~$G$, and let $H\preceq G$ mean that graph~$H$
can be obtained from graph~$G$ by contracting some of its edges.

By $C(H,k)$ we denote the chromatic polynomial of graph
$H$, that is, the number of proper colorings of the vertices
of graph~$H$ in at most~$k$ colors, and by $F(H,k)$ we denote
the flow polynomial of graph~$H$, that is the number of
flows modulo~$k$ having neither sources nor sinks and
equal to~$0$ on none of the edges.

\

2. RESULTS. We shall prove one theorem and deduce three of
its corollaries.

\

THEOREM. \emph{For every graph G}
\begin{equation}
C(G,k)=\frac{(k-1)^{m(G)}}{k^{m(G)-n(G)}}
\sum_{H\le G}\frac{F(H,k)}{(1-k)^{m(H)}}.
\label{eq1}\end{equation}

\

EXAMPLE 1. Let $G$ be a tree, then
$m(G)=n(G)-1$. If $H\le G$ then $F(H,k)=0$
with exception of the degenerate case
$m(H)=n(H)=0$ when $F(H,k)=1$. Hence,
$$ C(G,k)=k(k-1)^{m(G)}.$$

\

EXAMPLE 2. Let $G$ be a simple circle,
then $m(G)=n(G)$. If $H\le G$ then $F(H,k)$
is different from zero only in two
extreme cases: when $m(H)=m(G)$, that is
when $H=G$, and when $m(H)=0$. In the former
case $F(H,k)=k-1$, in the latter case
$F(H,k)=1$; hence
$$C(G,k)=(-1)^{m(G)}(k-1)+(k-1)^{m(G)}.$$

In \cite{ref1} the following representation
was obtained:
\begin{equation}
C(G,k)=\frac{(k-1)^{m(G)}}{k^{m(G)-n(G)}} \sum_{H\le
G}\frac{w(H,k)}{(1-k)^{m(H)}}. \label{eq2}\end{equation}
The
function~$w$ was  first defined as the sum of the values of a
certain weight function with summation over all (that is, non
necessary proper) colorings of graph~$G$ (formula 2.5 in
\cite{ref1}); later $w$ was defined via a recurrent relation
(Theorem VII in \cite{ref1}). An easy analysis of the proofs shows that
representations
\eqref{eq1} and \eqref{eq2}
termwise coincide, that is, always we have
\begin{equation}
F(H,k)=k^{m(H)-n(H)}w(H,k).
\label{eq3}\end{equation}
Thus the proposed theorem can be viewed as a
relationship between the special  function
introduced in \cite{ref1} and more traditional
notions in graph theory.

The definition of the function~$w$ implied by \eqref{eq3}
makes evident a number of its properties established in
\cite{ref1} to facilitate computation of the
function:
\begin{itemize}
\item $w(H,k)=0$ as long as $G$ contains an isthmus
(Theorems I and IV in \cite{ref1});
\item $w(H,k)=w(H_1,k)w(H_2,k)$ provided that $H$
consists of two parts without common vertices
(Theorem II in \cite{ref1});
\item $w(H,k)=w(H_1,k)w(H_2,k)/k$ provided that $H_1$
and $H_2$ have a single common vertex
(Theorem III in \cite{ref1});
\item $w(H',k)=w(H'',k)$ provided that $H'$
and $H''$ are homeomorphic
(Theorem V in \cite{ref1}).
\end{itemize}

The transition from \eqref{eq2} to \eqref{eq1} is
most interesting when $G$ is a planar graph. In this
case each subgraph~$H$ in \eqref{eq1} is also planar
and we can find its geometric dual graph~$H^*$. It is
easy to check that
\begin{equation}
F(H,k)=C(H^*,k)/k
\label{eq4}\end{equation}
and hence we have

\

COROLLARY 1. For every planar graph~$G$
\begin{equation}
C(G,k)=\frac{(k-1)^{m(G)}}{k^{m(G)-n(G)+1}}
\sum_{H\le G}\frac{C(H^*,k)}{(1-k)^{m(H)}}.
\label{eq5}\end{equation}

\

This result can be restated in a dual form:

\

COROLLARY 1. For every planar graph~$G$
\begin{equation}
C(G,k)=\frac{(k-1)^{m(G^*)}}{k^{n(G^*)-s}}
\sum_{L\le G^*}\frac{C(L,k)}{(1-k)^{m(L)}}.
\label{eq6}\end{equation}

\

(This result shows, in particular, how one can find the
chromatic polynomial of a connected planar graph~$G$
from its combinatory dual graph~$G^*$, although the
graph~$G$ itself isn't, in general,  determined
uniquely by~$G^*$.)

If $m(H)\ge 1$ then it is easy to see that
$C(H,k)\equiv 0 \pmod{k-1}$. This implies that,
in the case when $m(G)>1$, passing from \eqref{eq5}
to congruence modulo~$(k-1)^2$, we cam omit all
summands except the one corresponding to the
case~$H=G$. Thus we have

\

COROLLARY 2. \emph{If $m>1$ then}
\begin{equation}
C(G,k)\equiv (-1)^mC(G^*,k)\pmod{(k-1)^2}
\label{eq7}\end{equation}

\

Putting here $k=3$, we get

\

COROLLARY 3. If a planar graph~$G$ is different from
the full graph~$K_2$ and has exactly one
(up to renaming of colors) proper coloring of vertices
in three colors, then the graph~$G^*$ dual to graph~$G$
is also vertex colorable in three colors.

\

3.DEFINITIONS. Let $G$, $H$ be graphs, and let
$V(G)$, $V(H)$, $E(G)$, $E(H)$ be the corresponding
sets of vertices and edges. We say that graph~$H$
is an edge subgraph of graph~$G$ (and write $H\le G$)
if $E(H)\subseteq E(G)$ and $V(H)$ consists of those
and only those vertices of~$G$ that are incident to
edges from~$E(H)$. For the sake of validity of
formula \eqref{eq1} we admit the case when $E(H)$
(and hence $V(H)$) is the empty set.

The operation of contracting graph~$G$ by edges connecting two
adjacent vertices~$v'$ and~$v''$ consists in removing those (and
only those) edges and identifying vertices~$v'$ and~$v''$ into a
single vertex~$v$; thus, if vertices~$v'$ and~$v''$ had been
connected to a vertex~$w$ by paths of $l'$ and $l''$ edges
respectively, then the new vertex~$v$ is connected to $w$ by a path
of $l'+l''$ edges. We say that graph~$H$ is a contraction of
graph~$G$ (and write
$H\preceq G$) if $H$ can be obtained from $G$ by a number, possibly
zero, of edge contractions.

We take the ring~$R_k$ of residues modulo~$k$ as the
standard set of $k$ colors. By a vertex coloring of
graph~$H$ we mean any function defined on $V(H)$
with values from~$R_k$; a coloring is called proper if
the ends of each edge have different colors. By
$S(H,k)$ we denote the set of all  colorings in
$k$ colors, and by $S^+(H,k)$ we denote the set of all
proper  colorings in
$k$ colors. In this notation
$C(H,k)=|S^+(H,k)|$ is the cardinality of
$S^+(H,k)$. It is well-known (see, for example, \cite{ref2})
that for a fixed~$H$, the  function~$C(H,k)$ is a polynomial
of degree~$n(H)$ with integer coefficients.

In order to be able to introduce flows on a graph~$H$
we need to fix some orientation of all edges, which will be
done by denoting by $e'$ and $e''$ the beginning and the
end of an edge~$e$ respectively. When such an  orientation is
fixed, a flow modulo~$k$ on graph~$H$ is defined as any
function defined on $E(H)$ with values from~$R_k$.
(It is supposed that if we change the orientation of some
edges, we'll have to change the  sign of the flow on those
edges; all notions introduced below are invariant with
respect to such transformations.)
We say that a flow~$t$ is balanced if it has neither
sources nor sinks, that is if for every vertex~$v$
\begin{equation}
\sum_{e'=v}t(e)=\sum_{e''=v}t(e)
\label{eq8}\end{equation}
(the summation is performed in the ring~$R_k$,
that is, modulo~$k$). The degree of degeneracy $d(t)$
of a flow~$t$ is defined as the number of edges~$e$ such
that~$t(e)=0$; a flow~$t$ is called non-degenerate if
$d(t)=0$. By $T(H,k)$ we denote the set of all flows
modulo~$k$ on a graph~$H$, and by $T^=(H,k)$ we
denote the set of all balanced flows.
The flow polynomial~$F(H,k)$ equal to the number
of non-degenerate balanced flows was introduced
in \cite{ref3} (in \cite{ref3} it was denoted
$\phi(H,k)$; see also \cite[Section 14C]{ref4}).

If $H$ is a plane graph, that is, we have fixed
a mapping of it to a plane, then its geometric dual graph
$H^*$ is defined in the following way. Vertices of $H^*$
correspond to the areas on which graph~$H$ cuts the plane,
and the edges of $H^*$ correspond to the edges of graph
$H$: an edge~$e^*$ from $E(H^*)$ connects vertices
$v_1^*$ and $v_2^*$ from $V(H^*)$ if the edge~$e$ dual to
$e^*$ separates the areas corresponding to vertices
$e_1^*$ and $e_2^*$ (isthmuses of graph~$H$ correspond to
loops in graph~$H^*$).

\

4.PROOFS. Let $G$ be an arbitrary graph, $k$ be a positive
integer. We will use the shorthand $V=V(G)$, $n=n(G)$,
$S=S(G,k)$, and so on.

Let $\varepsilon=\cos(2\pi/k)+\mathrm{i}\sin(2\pi/k)$ be a primitive
root of unity of degree~$k$. Then for $r\in R_k$
\begin{equation*}
\sum_{t=0}^{k-1}\varepsilon^{rt}= \left\{\begin{array}{rcl}
k,&\mathrm{\ if\ }r=0,\\
0,&\mathrm{\ if\ }r\neq 0.
\end{array}\right.
\end{equation*}
From this we get that for $s\in S$
\begin{equation}
\prod_{e\in E}\left(k-\sum_{t=0}^{k-1}
\varepsilon^{(s(e')-s(e''))t}\right)= \left\{\begin{array}{rcl}
k^m,&\mathrm{\ if\ }s\in S^+,\\
0,&\mathrm{\ if\ }s\not\in S^+.
\end{array}\right.
\label{eq9}\end{equation}
Let
\begin{equation*}
\delta(t)=
\left\{\begin{array}{rcl}
k-1,&\mathrm{\ if\ }t=0,\\
-1,&\mathrm{\ if\ }t\neq 0.
\end{array}\right.
\end{equation*}
then by \eqref{eq9}
\begin{equation*}
k^mC(G,k)=
\sum_{s\in S^+}k^m=
\sum_{s\in S}\prod_{e\in E}\sum_{t=0}^{k-1}
   \delta(t)\varepsilon^{(s(e')-s(e''))t}.
\end{equation*}

Further we have:
\begin{equation*}
\prod_{e\in E}\sum_{t=0}^{k-1}
   \delta(t)\varepsilon^{(s(e')-s(e''))t}=
\sum_{t\in T}\prod_{e\in E}
   \delta(t(e))\varepsilon^{(s(e')-s(e''))t(e)},
\end{equation*}
\begin{equation*}
\prod_{e\in E}
   \delta(t(e))\varepsilon^{(s(e')-s(e''))t(e)}
=\prod_{e\in E}\delta(t(e))\times \prod_{e\in
E}\varepsilon^{s(e')t(e)}\times \prod_{e\in
E}\varepsilon^{-s(e'')t(e)},
\end{equation*}
\begin{equation*}
\prod_{e\in E}
   \delta(t(e))=(-1)^m(1-k)^d(t),
\end{equation*}
\begin{eqnarray*}
\prod_{e\in E}\varepsilon^{s(e')t(e)}
&=&\prod_{\parbox{8mm}{\scriptsize $e \in E$\\ $v\in V$\\$ e'=v$}}\varepsilon^{s(v)t(e)}\\
&=&\prod_{\ v\in V}\prod_{e'=v}\varepsilon^{s(v)t(e)}\\
&=&\prod_{\ v\in V}\varepsilon^{\sum_{e'=v}s(v)t(e)}\\
&=&\prod_{\ v\in V}\varepsilon^{\sum_{e'=v}t(e)\times s(v)}.
\end{eqnarray*}
Similarly,
\begin{eqnarray*}
\prod_{e\in E}\varepsilon^{-s(e'')t(e)} &=&\prod_{\ v\in
V}\varepsilon^{-\sum_{e''=v}t(e)\times s(v)},
\end{eqnarray*}
so that,
\begin{eqnarray}
k^mC(G,k)
&=&\sum_{s\in S,\ t\in T}(-1)^m(1-k)^{d(t)}
   \prod_{\ v\in V}\varepsilon^
   {\left(\sum_{e'=v}t(e)-\sum_{e''=v}t(e)\right)s(v)}\nonumber\\
&=&(-1)^m\sum_{t\in T}(1-k)^{d(t)}
    \prod_{\ v\in V}\sum_{s=0}^{k-1}
       \varepsilon^   {\left(\sum_{e'=v}t(e)-\sum_{e''=v}t(e)\right)s(v)}\nonumber\\
&=&(-1)^m\sum_{t\in T^=}(1-k)^{d(t)}k^n\nonumber\\
&=&(-1)^mk^n\sum_{t\in T^=}(1-k)^{d(t)}.
\label{eq10}
\end{eqnarray}

For every flow~$t$ from $T^=$ we define the
subgraph~$G_t$ as the graph obtained from $G$ by
removing those and only those edges~$e$ for which~$t(e)=0$;
clearly, $d(t)=n(G)-m(G_t)$. It is easy to see that
the restriction of a flow~$t$ on the graph~$G_t$ is a
balanced non-degenerate flow on~$G_t$, and, vice versa, for
every balanced non-degenerate flow~$t_H$ on any spanning subgraph
$H$ there exists a unique flow~$t$ such that $G_t=H$ and
$t_H$ is the restriction of $t$ on~$H$. Thus
\begin{equation*}
F(H,k)=\sum_{t\in T^=(G,k),\ G_t=H}1.
\end{equation*}

Continuing from \eqref{eq10}:
\begin{eqnarray*}
k^mC(G,k)
&=&(-1)^mk^n\sum_{t\in T^=}(1-k)^{d(t)}\\
&=&(-1)^mk^n\sum_{H\le G}\sum_{t\in T^=(G,k),\ G_t=H}(1-k)^{d(t)}\\
&=&(-1)^mk^n\sum_{H\le G}(1-k)^{-m(H)}\sum_{t\in T^=(G,k),\ G_t=H}1\\
&=&(-1)^mk^n\sum_{H\le G}\frac{F(H,k)}{(1-k)^{m(H)}}.
\end{eqnarray*}

The Theorem is proved.

Relation \eqref{eq4} is given in \cite{ref3} without proof
(see also \cite[Section 14C]{ref4}). For completeness we prove it now.

The addition operation of the ring~$R_k$ induces an addition in
$T(H,k)$; namely, let $t_1+t_2$ be such a flow that
$(t_1+t_2)(e)=t_1(e)+t_2(e)$ for $t_1,t_2\in T(H,k)$ and $e\in
E(H)$. Similarly, the multiplication operation of $R_k$ allows us to
multiply the elements of $T(H,k)$ by the elements of this ring:
$(rt)(e)=r\cdot t(e)$ for $t\in T(H,k)$, $r\in R_k$, and $e\in
E(H)$. Thus we can view $T(H,k)$ as a module over the ring~$R_k$. It
is easy to check that if $t,t_1,t_2\in T^=(H,k)$, then
$t_1+t_2,rt\in T^=(H,k)$; hence $T^=(H,k)$ is a submodule
of~$T(H,k)$.

If $H$ has an isthmus,  then $H^*$ has a loop, and thus
$F(H,k)=0=C(H^*,k)$; from now on we assume that $H$ has
no isthmus.

Let us imagine that graphs~$H$ and $H^*$ are drawn on a sphere. Let
$v^*$ be a vertex of graph~$H^*$ corresponding to a certain area
among areas on which graph~$H$ divides the whole sphere. Let
$e_1,\dots,e_q$ be the edges bounding this area. For the definition
of flows these edges were somehow oriented, so now we can speak of
the edges~$e_1,\dots,e_q$ as oriented clock-wise and
counterclock-wise (assuming that the ``center of the clock'' is at
the vertex~$v^*$). Let us define a flow~$t_{v^*}$ as the flow equal
to~$+1$ on clock-wise oriented edges, $-1$ on edges oriented in the
opposite direction, and equal to~$0$ on the remaining edges (that
is, different from ~$e_1,\dots,e_q$). Clearly, the flow~$t_{v^*}$ is
balanced.

To a given coloring~$s^*$ from $S(H^*,k)$, we associate
the balanced flow:
\begin{equation}
t_{s^*}=\sum_{v^*\in V(H^*)}s^*(v^*)t_{v^*}.
\label{eq11}\end{equation}
Let us show that for each balanced flow~$t$ on $H$ there
exist exactly $k$ colorings~$s^*$ such that
\begin{equation}
t=t_{s^*}.
\label{eq12}\end{equation}

First we prove that the number of such colorings cannot
be greater than~$k$. To this end we fix a vertex~$v_0$
from $V(H^*)$ and show that a coloring~$s^*$ satisfying
condition \eqref{eq12} can be uniquely determined
by its values~$s^*(v_0^*)$. Because of the connectivity
of the graph~$H^*$ it is sufficient to show that
the value ~$s^*(v_1^*)$ is uniquely determined where
$v_1^*$ is a vertex adjacent to~$v_0^*$. Let $l$ be
the edge dual to the edge connecting $v_0^*$ and~$v_1^*$.
According to \eqref{eq12} and \eqref{eq11}
\begin{equation}
t(l)=t_{s^*}(l)=\pm (s^*(v_0^*)-s^*(v_1^*))
\label{eq13}
\end{equation}
(The sign depends on the orientation of the edge~$l$),
and this relation allows us to determine
$s^*(v_1^*)$ from  $s^*(v_0^*)$ and~$t$.)

Because $s^*(v_0^*)$ can assume at most
$k$ values, the number of colorings satisfying \eqref{eq12}
cannot be greater than~$k$. Let us now show that indeed
all $k$ cases can be implemented.

Let us fix a spanning tree~$D^*$ of graph~$H^*$. Let us
take for the value of $s^*(v_0^*)$ an arbitrary element
of the ring~$R_k$ and define values of $s^*$ on other
vertices of graph~$H^*$ according to the above described
procedure using relation \eqref{eq13} only for edges
dual to the edges of the tree~$D^*$. Let us show
that the resulting coloring will satisfy the equality
$t(l)=t_{s^*}(l)$ for the other edges as well. These
edges form a tree~$W$. By construction of $s^*$
the flow~$t-t_{s^*}$ is equal to zero outside~$W$,
and thus its restriction to $W$ is balanced as well.
But the only balanced flow on a forest is the flow
identically equal to zero.

In order to complete the proof of equality \eqref{eq4}
it remains to note that a coloring is proper if and only
if the flow corresponding to it is non-degenerate.

In order to pass from \eqref{eq5} to \eqref{eq6} it
suffices to note that by definition
$m(G)=m(G^*)$, by Euler's Theorem
$n(G)-m(G)+n(G^*)=1+s$, and there is a natural
one-to-one correspondence between the sets
$\{H^*|H\le G\}$ and $\{L|L\preceq G^*\}$:
namely, $H^*$ is obtained from $G^*$ by
contracting edges dual to edges from~$V(G)\setminus V(H)$.

\

\hfill \parbox{60mm}{Received by the editorial board\\
on April  1, 1977}

\end{document}